\documentclass[preprint,12pt]{elsarticle}

\usepackage{hyperref}
\usepackage{amstext}
\usepackage{graphicx}
\usepackage{epstopdf}
\usepackage{amssymb}
\usepackage{mathrsfs}
\usepackage{amsmath}
\usepackage{amsthm}
\usepackage{amsfonts}
\usepackage{verbatim}
\usepackage{color}
\usepackage{natbib}
\usepackage{multirow}

\newcommand{\bbR}{\mathbb{R}}

\newcommand{\pd}{\partial}

\journal{Journal of Computational Physics}

\begin{document}

\begin{frontmatter}

\title{Interface Formulation and High Order Numerical Solutions of PDEs with Low Regularity} 

\author{Y. C. Zhou\corref{mycorrespondingauthor}}
\address{Department of Mathematics, Colorado State University, Fort Collins, CO 80523, USA}
\cortext[mycorrespondingauthor]{Corresponding author}
\ead{yzhou@math.colostate.edu}

\author{Varun Gupta}
\address{Pacific Northwest National Laboratory, Richland, WA 99352, USA}




\begin{abstract}
Linear elastic fracture mechanics admit analytic solutions that have low regularity 
at crack tips. Current numerical methods for partial differential equations (PDEs) of this type 
suffer from the constraint of such low regularity, and fail to deliver optimal high order rate of convergence. 
We approach the problem by (i) choosing an artificial interface to enclose the center of the low regularity; and (ii)
representing the solution in the interior of artificial interface as unknown linear combination of known modes of
low regular solutions. This gives rise to an interface formulation of the original PDE, and the linear combination 
are represented the interface conditions. By enforcing the smooth component of numerical solution in the interior 
domain to be approximately zero, a least square problem is obtained for the unknown coefficients. The solution of 
this least square problem will provide approximate interface conditions for the numerical solution of the PDE in the exterior 
domain. The potential of our interface formulation is favorably demonstrated by numerical experiments
on 1-D and 2-D Poisson equations with low regular solutions. High order numerical solutions of
unknown coefficients and PDEs are obtained. This proves the potential of the proposed interface formulation as 
the theoretical basis for solving linear elastic fracture mechanics problems. 
We indicate the relations between our interface formulation and domain decomposition methods as well as a regularization 
strategy for the Poisson-Boltzmann equation with singular charge density.
\end{abstract}

\begin{keyword}
Solution of low regularity \sep approximated interface condition \sep least square solution \sep high order
\end{keyword}

\end{frontmatter}


\section{Introduction}
This study is motivated by the challenges in the construction of high order numerical methods for solving linear elastic 
fracture mechanics problems. A classical example is the propagation of a crack in the toughness-dominated problems. 
For certain geometric configurations and boundary conditions, analytical solutions of the displacement fields near the crack front 
(referred to as crack tip in 2-D) can be obtained \cite{AndersonFractureMechanics}. These solutions contain 
terms of the form $r^p$ with $\frac{1}{2} \le p \le \frac{3}{2}$, and thus the stress field is singular at the front where $r=0$.
These fracture problems are the prototypes of the so-called L-shape elliptic problems with similar 
solution singularities at the boundary \cite{SzaboBabuska}. Standard finite element or boundary integral methods
with uniform mesh refinement yield sub-optimal 
convergence on this type of problems, and thus various modifications of numerical methods have been developed to 
recover the formal order of asymptotic convergence of the numerical solutions, including 
local adaptive techniques \cite{OdenJ1998a,ElliotisM2002a,NochettoR2011a,BadiaS2013a,MitchellW2015a}.

Extended finite element methods (XFEM) (also called generalized finite element methods (GEFM)) 
are the most promising numerical techniques that can capture the singular nature of the stress field in the 
numerical solutions \cite{BabuBanOsbGFEM, StrouboulisT2001a,FriesT2010a,xfem_crack2D, ABO_CS}.
XFEM enhances the solution space of the finite element
methods by adding enrichment functions to account for the discontinuous displacement across the crack surface and
the singular stress field at the crack front. In particular, the known analytical solution of the displacement 
field with low regularity is adopted as the enrichment function at the nodes in the vicinity of
the crack front \cite{bechet:improvedXFEM:2005,MoesN1999a,ABO_CS,Gupta2015355}. A salient
feature of this enrichment scheme is that one can use a fixed mesh and adaptively add enrichment functions
only around the crack following the propagation of crack on a fixed finite element mesh. 
Algebraically, this means one can have a major trunk of the stiffness matrix corresponding to the
basic degrees of freedom unchanged during the growth of crack, and only requires updating a small
portion of the matrix constructed by hierarchically adding enrichment degrees of freedom. Optimal convergence of linear XFEM has been observed 
for crack problems \cite{laborde:HighOrderCrack:2005,
bechet:improvedXFEM:2005,Renard_Chahine_optimal_convergence_xfem_2011,Tarancon2009,
Gupta2016481,Gupta2015355}, and its mathematical analysis has been well 
established \cite{BabuBanOsbGFEM, StrouboulisT2001a}. It is natural to construct and apply high order XFEM to crack problems, which
nevertheless, failed to deliver consistent optimal rate of convergence through 
numerical experiments \cite{laborde:HighOrderCrack:2005,ZamaniA2010a,ChahineE2011a,Gupta2016481}.

The use of fixed mesh and known modes of low regular solution at the crack front in XFEM
motivates us to investigate whether it is possible to formulate interface problems
for such crack problems. That way various numerical methods 
can be deployed to the resultant interface problems to produce high order numerical solutions
that are otherwise hard to obtain using high order extensions of XFEM. In fact, we could generalize this idea to other 
elliptic problems with known modes of lower regularity in their solutions. Note that this does not mean that the 
solution is completely known at the position where low regularity is observed. In general these 
known modes of low regular solutions are found with special geometric configurations and 
boundary conditions, and their linear combination is to be determined for the specific boundary conditions 
of a given problem. The superposition rule also allows the existence of other regular solution components in 
addition to the low regular solutions at any point in the domain. As described below, we are going to define the solution 
as a linear combination of known modes of low regular solutions 
and modes of smooth, regular solutions. We will subtract this solution representation in a chosen region around the 
tip to create an interface problem whose jump conditions are given by the constructed solution. The 
discretization of the resulting interface problem contains the unknown coefficients of linear combination. 
We present a least-square procedure to solve for these coefficients, which will be then supplied to  
solve the interface problem. Our construction is independent of the dimensionality of the problem
and of the numerical methods for solving the associated interface problems. Below we will present the detail 
of the interface formulation, its numerical analysis, and computational experiments on 1-D and 2-D problems. 

\section{Interface Formulation and its Solution} \label{sect:model}
We consider an open domain $\Omega \in \bbR^n$ with piecewise Lipschitz boundary $\pd \Omega$ and a Poisson equation
with Dirichlet boundary condition
\begin{equation} \label{eqn:model_problem}
\begin{split}
\nabla \cdot (a(x) \nabla u) & = f(x), \quad x \in \Omega  \\
\quad u(x) & = g(x), \quad x \in \partial \Omega
\end{split}
\end{equation} 
We are interested in Eq.(\ref{eqn:model_problem}) with a source function $f(x)$ or a boundary value $g(x)$ 
whose singularity at $x_0 \in \bar{\Omega}$ leads to the low regularity of the solution $u(x)$ at $x_0$. 

\begin{figure}[!ht] 
\centering 
\def\svgwidth{0.80\textwidth}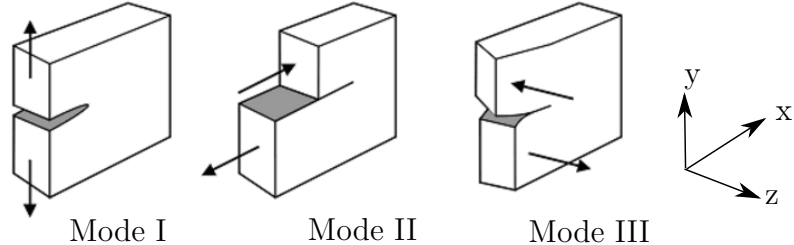
\caption{ Three basic modes of fracture: the opening (Mode I), the in-plane shear (Mode II)
and the out-of-plane shear (Mode III).}
\label{fig:fractureModes}
\end{figure}
A concrete example of Eq.(\ref{eqn:model_problem}) is the linear elastic fracture mechanics, for which 
one can identify three independent movements of the disjoint  
upper and lower faces of a crack during the propagation, generating three fundamental modes of linear
elastic fracture as shown in Figure \ref{fig:fractureModes}
\cite{AndersonFractureMechanics}. Furthermore, one can take advantage of the
symmetry of all three crack deformation modes to reduce their mathematical description to 
two-dimensional linear elastic problems. All three modes have solution component with $\sqrt{r}$
term, leading to singularity at $r=0$. As an example, the $x-$component of Mode I displacement in the local crack front coordinate system is given by,
\begin{equation}
\label{eq_ez_size:mode_1_expansion}
\left(\mathbf{u}_I(r,\theta)\right)_x = \sqrt{r} \left\lbrace
\left(\kappa-{1\over2} \right)\cos{\theta \over 2} - {1 \over 2}\cos{3\theta
\over 2}\right\rbrace 
\end{equation}
%
where $r$ and $\theta$ are polar coordinates at the crack front,  $-\pi <
\theta < \pi$, $\kappa$ is the material constant ($3-4\nu$), and $\nu$ is the
Poisson's ratio. Complete elasticity solution in the neighborhood 
of a crack can be found in \cite{AndersonFractureMechanics}.

Without loss of generality we assume that for Eq.(\ref{eqn:model_problem})
the solution of low regularity at $x_0$ is given by a linear combination of finite number of modes, and the full solution 
near $x_0$ can be given by this solution of low regularity adding up to some regular components, in the following general form
\begin{equation}
u_s(x) = \sum_{i=1}^{N_l} k^l_i u^l_i(x) + \sum_{j=1}^{N_r} k^r_j u^r_j(x)
\end{equation}
where $u^l_i(x)$ is the $i^{th}$ mode of the low regular solution, and $u^r_j(x)$ is the $j^{th}$ component of regular solution. The
set $\{ u^l_i(x) \}$ must be linearly independent, so is the set $\{ u^r_j(x) \}$. 
We will be concentrated in the treatment of the low regular solutions, whose known modes are typically obtained by 
solving Eq.(\ref{eqn:model_problem}) with specialized source term and boundary conditions. 
We will discuss the construction of the regular components later at the end of Section \ref{sect:numerical}.  
The regular components can be chosen from suitable polynomial space, 
such as those that form an orthogonal basis of $L^2(\Omega_1)$ where $\Omega_1 \subset \Omega$ 
is a small domain in which $x_0$ is located. For the purpose of analysis it is not necessary to distinguish $u^l_i$ and $u^r_j$, 
so a general expansion formulation of the solution around $x_0$ 
\begin{equation} \label{eqn:us_all}
u_s(x) = \sum_{i=1}^{N_l+N_r} k_i u^s_i(x)
\end{equation}
is adopted in the discussion below.

We will choose a subdomain $\Omega_1 \subset \Omega$ to enclose $x_0$. We will subtract $u_s(x)$ from the solution $u(x)$ in $\Omega_1$, 
while in the $\Omega_2 = \Omega \setminus \Omega_1$ the original solution $u(x)$ is untouched. We call this subtracted solution 
as $\tilde{u}(x)$, which will solve the following interface problem
\begin{equation} \label{eqn:interface_problem}
\begin{split}
\nabla \cdot (a(x) \nabla \tilde{u}) = f(x) - \sum_{i=1}^{N_l+N_R} k_i \nabla \cdot(a(x) \nabla u^s_i(x)), \quad x \in \Omega_1, \\  
\nabla \cdot (a(x) \nabla \tilde{u}) = f(x),  \quad x \in \Omega_2, \\  
\quad \tilde{u}(x) = g(x), \quad x \in \partial \Omega \\
\left [ \tilde{u} \right ] = \sum_{i=1}^{N_l+N_r} k_i u^s_i(x), \quad x \in \Gamma, \\
\left [ a(x) \nabla \tilde{u} \right ] = a(x) \sum_{i=1}^{N_l+N_r} k_i \nabla u^s_i(x), \quad x \in \Gamma,
\end{split}
\end{equation} 
where $[\cdot ]$ denotes the jump at the interface $\Gamma$ of the enclosed quantity from $\Omega_1$ to $\Omega_2$. {\it Our goal
is to solve $\tilde{u}(x)$ in $\Omega$ and the expansion coefficients $k_i$}.

The past decades have witnessed the development of many numerical methods for approximating the
elliptic interface problems in the form of Eq.(\ref{eqn:interface_problem}). 
The following general linear algebraic system can be obtained after the application
of one of these interface methods to Eq.(\ref{eqn:interface_problem}):
\begin{equation} \label{eqn:linear_system}
\mathcal{S} \tilde{u} = F + E k \Leftrightarrow
\left [  
\begin{array}{cc}
A & B \\
C & D 
\end{array}
\right ] 
\cdot
\left [ 
\begin{array}{c}
\tilde{u}_1 \\
\tilde{u}_2 
\end{array}
\right ] 
= 
\left [ 
\begin{array}{c}
F_1 \\
F_2 
\end{array}
\right ] 
+
\sum_{i=1}^{N_l + N_r}  k_i
\left [ 
\begin{array}{c}
E_{i,1}  \\
E_{i,2}  
\end{array}
\right ], 
\end{equation} 
where $\tilde{u}_1, \tilde{u}_2$ are the numerical solutions of $\tilde{u}(x)$ in subdomains $\Omega_1, \Omega_2$, respectively, and $A,B,C,D,F_1,F_2$
are the corresponding blocks in the stiffness matrix $\mathcal{S}$ and the loading vector $F$. In particular, $E_{i,1}$ and $E_{i,2}$ are 
the blocks of $E_i$, the contribution from the $i^{th}$ known mode $u^s_i(x)$ to the loading vector. The coefficients $k_i$ are unknown, but 
we expect the solution $\tilde{u}_1=0$ by checking the interface problem Eq.(\ref{eqn:interface_problem}). Therefore {\it our goal at this point is 
to find $\{ k_i\}$ and $\tilde{u}_2$ such that $\tilde{u}_1$ is (approximately) zero}. 

We could move the unknown contributions from the interface 
condition in the linear system (\ref{eqn:linear_system}) to the left hand size to obtain
\begin{equation} \label{eqn:linear_system_2}
\left [ 
\begin{array}{ccc}
A & B & E_1 \\
C & D & E_2 
\end{array}
\right ] 
\cdot
\left [ 
\begin{array}{c}
\tilde{u}_1 \\
\tilde{u}_2  \\
k
\end{array}
\right ] 
= 
\left [ 
\begin{array}{c}
F_1 \\
F_2 
\end{array}
\right ],
\end{equation} 
where $E_1,E_2$ and $k$ collect the associated components. Arguing that $\tilde{u}_1=0$ one can remove those degrees of freedom from 
Eq.(\ref{eqn:linear_system_2}) but the remaining system is still un-determined. Indeed it is always overdetermined practically 
if the degrees of freedom in $\Omega_1$ is larger than the number of known modes $N_l + N_r$. Taking the opening mode of the crack
problem as an example, $N_l = 1$ but one can easily deploy tens of nodes in a reasonably chosen domain $\Omega_1$ that encloses the tip $x_0$.

We now present two different approaches to solve $k$ and $\tilde{u}_2$ from Eq.(\ref{eqn:linear_system_2}). 
The first approach solving for $k$ and $\tilde{u}_2$ follows the solution of $\tilde{u}_1$ 
from Eq.(\ref{eqn:linear_system}) using its Schur complementary: 
\begin{equation} \label{eqn:schur_a}
\left [  
\begin{array}{cc}
A-BD^{-1}C & 0 \\
C & D 
\end{array}
\right ] 
\cdot
\left [ 
\begin{array}{c}
\tilde{u}_1 \\
\tilde{u}_2 
\end{array}
\right ] 
= 
\left [ 
\begin{array}{c}
F_1 - BD^{-1} F_2\\
F_2 
\end{array}
\right ] 
+
\sum_{i=1}^{N_l + N_r}  k_i
\left [ 
\begin{array}{c}
E_{i,1} - BD^{-1} E_{i,2} \\
E_{i,2}  
\end{array}
\right ], 
\end{equation}
i.e., 
\begin{equation} \label{eqn:schur_b}
(A - BD^{-1} C) \tilde{u}_1 = (F_1 - BD^{-1} F_2) + (E_1 - BD^{-1}E_2) k.
\end{equation}
In order to have a zero solution of $\tilde{u}_1$ we shall have zero right-hand side in Eq.(\ref{eqn:schur_b}), and that leads to a
least-square solution of $k$, which can be formally written as 
\begin{equation} \label{eqn:optimal_k_b}
k = -((E_1 - BD^{-1} E_2)^T (E_1 - BD^{-1} E_2))^{-1} (E_1 - BD^{-1} E_2)^T (F_1 - BD^{-1} F_2).
\end{equation}
With this $k$ we can solve $\tilde{u}_2$ from the following subsystem 
of Eq.(\ref{eqn:linear_system}):
\begin{equation} \label{eqn:u2_only}
D \tilde{u}_2 = F_2 + \sum_{i=1}^{N_l + N_r} k_i E_{i,2} = F_2 + E_2 k.
\end{equation}

In Eqs.(\ref{eqn:schur_a}-\ref{eqn:optimal_k_b}) $D^{-1}F_2, 
D^{-1}E_2$ are respectively the solutions of $Dx = F_2, Dy = E_2$, and they can be obtained using the same efficient iterative 
algorithm for solving $\tilde{u}_2$ from Eq.(\ref{eqn:u2_only}). We will add back $u_s(x)$ as defined in Eq.(\ref{eqn:us_all}) with
numerically solved $k_i$ to the interior domain $\Omega_1$ to get $\tilde{u}_1$.

An alternative strategy involves an iterative scheme for solving $k$ but it turns out that the mapping 
from $k^{j}$ to $k^{j+1}$ has a fixed point, which can be solved directly to given an optimal $k$, and then
the solution of $\tilde{u}_2$ follows. Choose some initial $k$ ($k^0 = 0$, for instance) and assume $k^j$ is 
known after $j$ steps of iteration, we solve $\tilde{u}$ from Eq.(\ref{eqn:linear_system}) as
\begin{equation} \label{eqn:linear_system_it_a}
\left [ 
\begin{array}{c}
\tilde{u}_1 \\
\tilde{u}_2 
\end{array}
\right ] 
= 
\mathcal{S}^{-1} 
(F 
+ E k^j),
\end{equation} 
where we suppress the block notions of $F$ and $E$. The $\tilde{u}_1$ block of $\tilde{u}$ is not necessarily zero, 
but can be made zero by multiplying from left with a proper matrix, i.e., 
\begin{equation} \label{eqn:linear_system_it_b}
\left [ 
\begin{array}{cc}
0 & 0 \\
0 & I 
\end{array}
\right ]. 
\left [ 
\begin{array}{c}
\tilde{u}_1 \\
\tilde{u}_2 
\end{array}
\right ] 
=
\left [ 
\begin{array}{c}
0 \\
\tilde{u}_2 
\end{array}
\right ] 
= 
\left [ 
\begin{array}{cc}
0 & 0 \\
0 & I 
\end{array}
\right ] \cdot 
\mathcal{S}^{-1} 
(F 
+ E k^j).
\end{equation}
We then multiply the original stiffness matrix $\mathcal{S}$ to this hacked solution $(0, \tilde{u}_2)^T$, and equalize the
product to $F + E k^{j+1}$, i.e.,  
\begin{align}
\mathcal{S} \cdot 
\left [ 
\begin{array}{c}
0 \\
\tilde{u}_2 
\end{array}
\right ]  
& = 
\mathcal{S} \cdot 
\left [ 
\begin{array}{cc}
0 & 0 \\
0 & I 
\end{array}
\right ] \cdot 
\mathcal{S}^{-1} 
(F 
+ E k^j) \nonumber \\
& = 
\left [ 
\begin{array}{cc}
0 & B \\
0 & D 
\end{array}
\right ] \cdot 
\mathcal{S}^{-1} 
(F 
+ E k^j)  \nonumber \\
& = F + E k^{j+1}.
\label{eqn:linear_system_it_c}
\end{align}
Notice that $\mathcal{S}^{-1} F$ is nothing but the solution of the linear system $\mathcal{S} x = F$, we denote $\mathcal{S}^{-1}F$ by $S^F$. 
We also collect the solutions of $\mathcal{S} x = E_i$ for all $i$ and denote it by $S^E$. With these notions we can write 
Eq.(\ref{eqn:linear_system_it_c}) as
\begin{equation*} 
F + E k^{j+1} = 
\left [ 
\begin{array}{c}
B S^F_1 \\
D S^F_2 \\
\end{array}
\right ] 
+  
\left [ 
\begin{array}{c}
B S^E_1 \\
D S^E_2 \\
\end{array}
\right ] k^j \equiv P + Q k^j, 
\end{equation*}
i.e.,
\begin{equation} \label{eqn:linear_system_it_d}
E k^{j+1} = P - F + Q k^j.
\end{equation}
This is an over-determined linear system for $k^{j+1}$, which can be solved in the least-square sense. Notice that one only needs to solve for
$S^F$ and $S^E$ once, before they are used repeatedly in the iterations for updating $k^{j+1}$. Interestingly, we can get around of this iteration by
directly calculating the fixed point of the mapping indicated by Eq.(\ref{eqn:linear_system_it_d}) using a stable least-square algorithm. This
solution can be formally written as 
\begin{equation}  \label{eqn:optimal_k_a}
k = ((E-Q)^T(E-Q))^{-1} ((E-Q)^T (P-F)).
\end{equation}
With this $k$ we can solve $\tilde{u}_2$ from equation \eqref{eqn:u2_only}. Eq.(\ref{eqn:optimal_k_b}) and Eq.(\ref{eqn:optimal_k_a}) 
provide similar solutions for $k$. In the computation below we adopt Eq.(\ref{eqn:optimal_k_b}) since the linear system it involves 
has a smaller dimensionality.

\section{Numerical Experiments} \label{sect:numerical}

We will test our interface formulation and solution techniques on 1-D and 2-D problems. Their extension to 3-D problems is conceptually straightforward.
Our interface formulation can be solved by using any interface method, and here we choose the Matched Interface and Boundary (MIB) method since
it comes with a convenient construction of high order approximations of interface conditions for attaining high order numerical solutions 
in the entire domain. Details of implementing MIB method are not relevant to the topics in this article, and interested readers are 
referred to \cite{ZhouY2006a,ZhouY2006b,YuS2007a,ZhouY2008a,ZhaoS2010a,WangB2015a} for its basic constructions and promising applications to a 
wide range of interface problems arising in electromagnetics, biomolecular electrostatics, and elasticity.

The first set of tests are on 1-D Poisson equation 
\begin{equation} \label{eqn:1d_poisson}
u_{xx} = f(x) 
\end{equation}
defined in the interval $[0,\pi]$. We choose an analytical solution
$u(x) = 2 x^{1/2}$ to compute $f(x)$ in the domain and the boundary conditions at two ends $x=0,\pi$. The solution has a single mode of low 
regular solution $u_s(x) = u^s_1(x) = x^{1/2}$ near $x=0$. We will use this $u_s(x)$ to construct an interface formulation of the form (\ref{eqn:interface_problem}),
and we choose $x_c = 0.5$ as the interface. This choice of $u_s(x)$ also sets the analytical solution of $k=2$. 
MIB methods of the second and fourth order are used for solving the interface formulation, along with the standard second order central 
finite difference method for the original Poisson equation (\ref{eqn:1d_poisson}). Rates of convergence are computed for all these three types of 
methods based on the numerical errors of ${u}$ evaluated in $L_{\infty}$ norm. The analytical solution of $k=2$
allows us to compare the error in numerical solutions of $k$ from the two MIB methods on the interface formulation of the Poisson equation.

\begin{table}[!ht]
\begin{center}
\begin{tabular}{rrrrrrrrr} \hline
\multirow{2}{*}{N} & \multicolumn{3}{c}{IF-MIB2}  & \multicolumn{3}{c}{IF-MIB4}  & \multicolumn{2}{c}{FD2}  \\ \cline{2-4} \cline{5-7} \cline{8-9}
         & $\| e \|_{\infty}$ & ROC & $|k-2|$  & $\| e \|_{\infty}$ & ROC  & $|k-2|$  & $\| e \|_{\infty}$ & ROC \\ \hline 
  20    & 1.00E-2  &       & 1.44E-2 & 2.63E-3  &        & 3.76E-3 & 2.61E-1  &   \\ 
  40    & 1.14E-3  & 3.13  & 1.48E-3 & 1.60E-4  & 4.04   & 2.28E-4 & 1.90E-1  & 0.46  \\ 
  80    & 7.37E-5  & 3.95  & 8.97E-5 & 7.35E-6  & 4.44   & 1.35E-6  & 1.38E-1  & 0.46  \\ 
  160   & 1.62E-5  & 2.19  & 1.77E-5 & 4.54E-7  & 4.02   & 6.42E-7 & 9.93E-2  & 0.47  \\ 
  320   & 4.17E-6  & 1.96  & 2.18E-6 & 5.03E-9  & 6.42   & 5.51E-9 & 7.09E-2  & 0.49  \\  
  640   & 1.06E-6  & 1.98  & 2.84E-7 & 2.2E-10  & 4.50  & 2.3E-10 & 5.05E-2  & 0.49  \\  \hline
\end{tabular}
\caption{Numerical error of ${u}$, its rate of convergence (ROC), and absolute error of $k$ 
for second-order and forth-order MIB methods with the interface formulation (IF-MIB2, IF-MIB4) and for 
second-order central finite difference method (FD2) on the 1-D Poisson equation with $u_s(x)=x^{1/2}$. 
$N$ is the number of uniform grid divisions.}
\label{tab:1d_poisson_error_u}
\end{center}
\end{table}

\begin{figure}[!ht]
\begin{center}
\includegraphics[width=6.5cm]{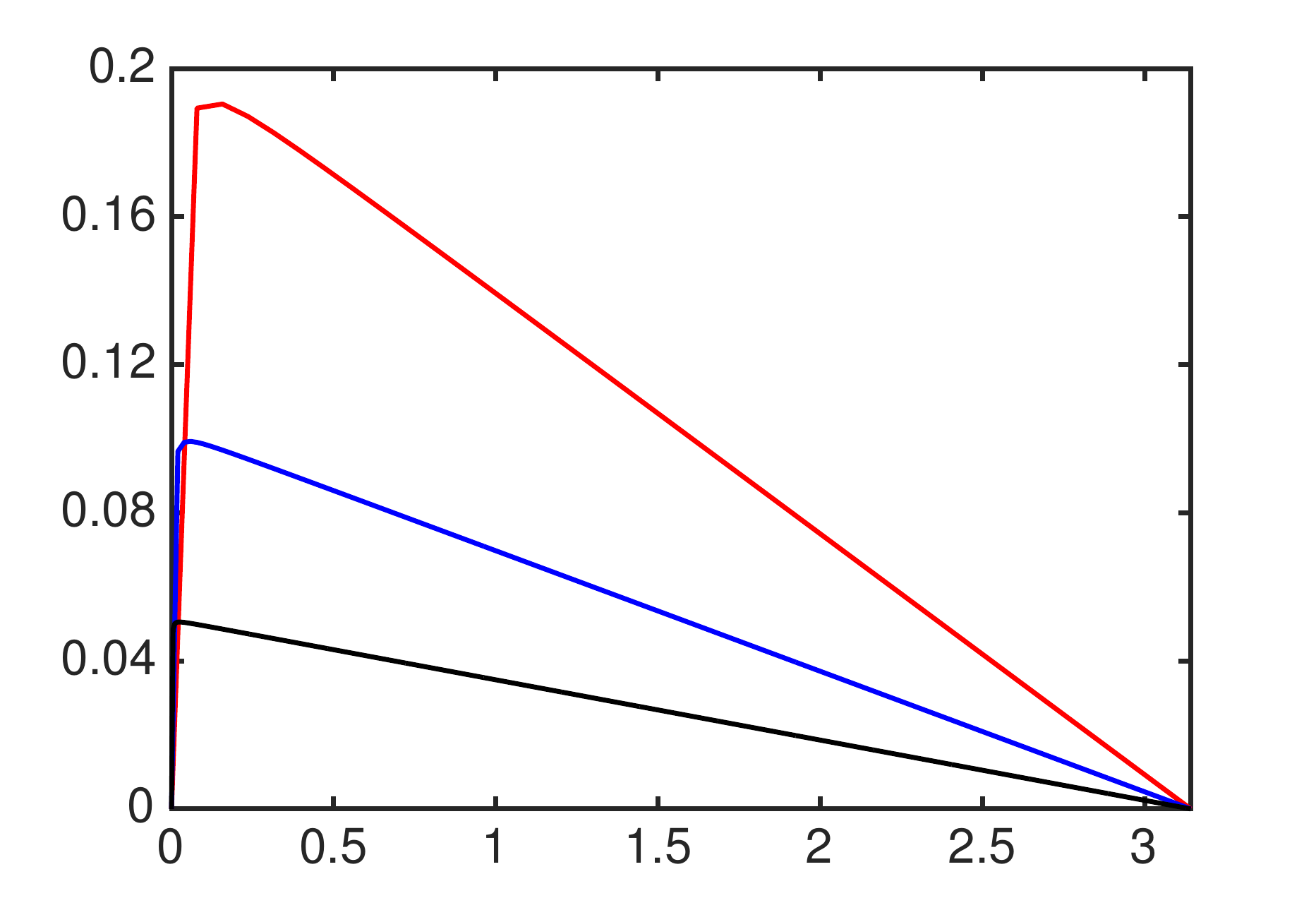}
\includegraphics[width=6.5cm]{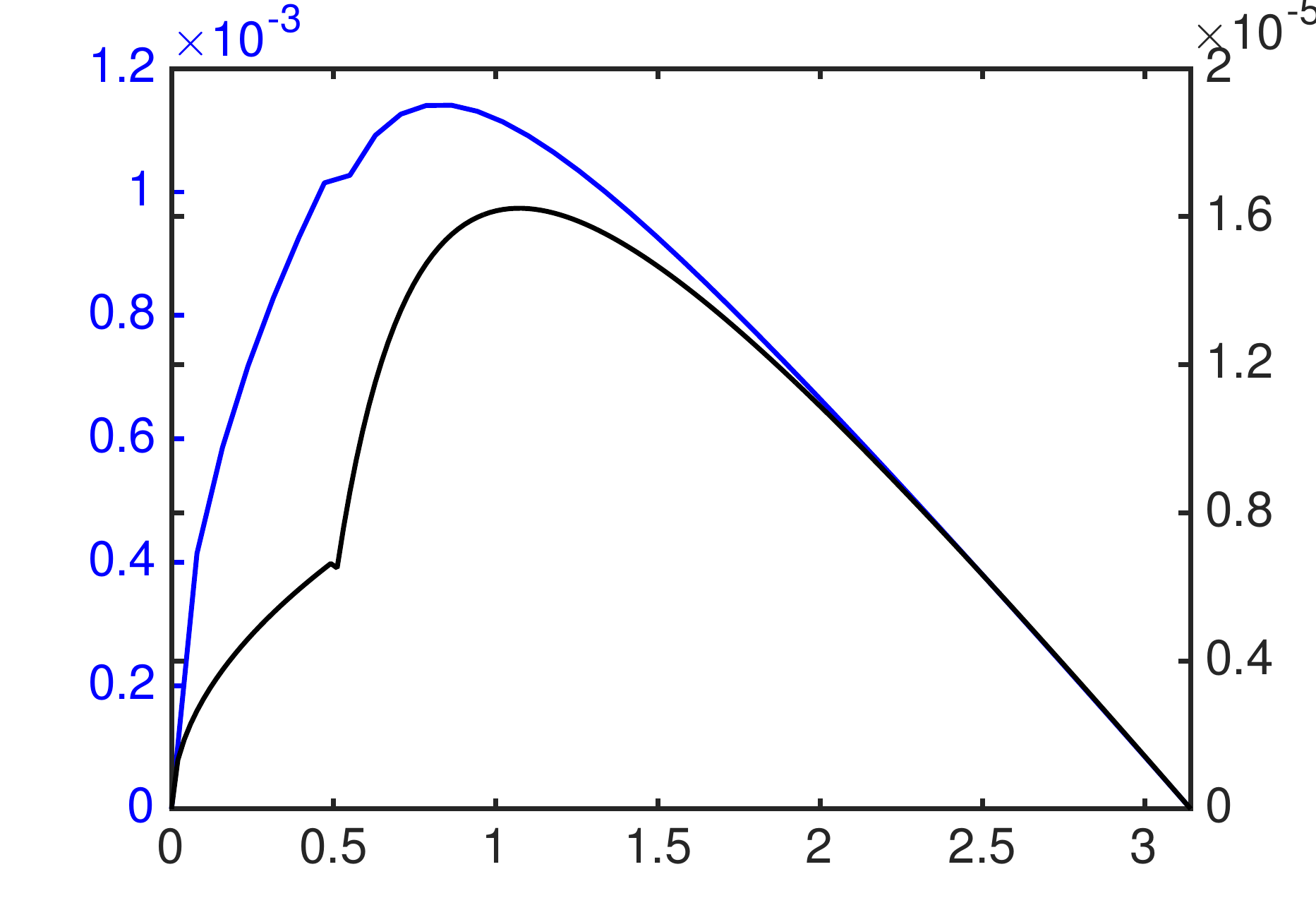}
\caption{Left: Errors for FD2 solutions at $N=40$(red), $80$(blue) and $160$(black). Right: Errors for
IF-MIB2 solutions at $N=40$(blue, left $y$-axis) and $N=160$ (black, right $y$-axis).}
\label{fig:error_1d_poisson}
\end{center}
\end{figure}

The $L_{\infty}$ errors collected in Table \ref{tab:1d_poisson_error_u} demonstrate that high order solution can be obtained from 
the interface formulation using the designed solution technique along with a suitable interface method. Not surprisingly the second-order 
central difference fails to deliver a second-order convergence for this Poisson equation due to the low regularity of $u(x)$ at $x=0$.
Notice that our interface formulation also involves an unknown coefficient $k$, whose absolute errors are also collected in 
Table \ref{tab:1d_poisson_error_u}. These data show that the error in $k$ is compatible with the error in ${u}$. Indeed,
the analytical relation in Eq.(\ref{eqn:schur_b}) suggest that for the exact value of $k$ it holds that
\begin{equation} \label{eqn:schur_analysis_a}
(E_1 - BD^{-1}E_2) k =  (F_1 - BD^{-1} F_2) - (A - BD^{-1} C) \tilde{u}_1. 
\end{equation}
By neglecting the term involving $\tilde{u}_1$ whose exact solution is zero, we are introducing to the right-hand size of
this least square problem a perturbation of $O(h^p)$ where $h$ is the mesh size and $p=2,4$ for IF-MIB2 and IF-MIB4, respectively. 
The stability condition of the least square problem will then assure the error in $k$ caused by this perturbation is also of $O(h^2)$. 
Fig.(\ref{fig:error_1d_poisson}) 
shows the concentration of the numerical error close to $x=0$ for second order central finite difference method, along with the
error distribution for the interface formulation solved with the second order MIB method. The error distributions for interface
solutions highlight that for this 1-D problem the interface conditions are very well approximated, and the maximum errors appear inside 
of the exterior domain.

We proceed to test the interface formulation with multiple modes of low regular solutions. We choose analytical solution of 
Eq.(\ref{eqn:1d_poisson}) to be $u(x) = 2 x^{1/2} + 3x^{3/2}$, and correspondingly we define $u^s_1 = x^{1/2}, u^s_2 = x^{3/2}$.
The numerical error of ${u}$ and rate of convergence are similar to what is shown in Table \ref{tab:1d_poisson_error_u},
so we only present the errors for the coefficients of two modes, as summarized in Table \ref{tab:1d_poisson_error_k_2}. The first
mode, $x^{1/2}$, dominates the other mode $x^{3/2}$ for small $x$, meaning the error in $k_1$ will contribute more to the
error in the numerical solution of $\tilde{u_1}$. Arguing inversely, an approximate solution of $k_1,k_2$ from the 
least-square procedure will lead to more accurate solution of $k_1$, as indicated in Table \ref{tab:1d_poisson_error_k_2}.
\begin{table}[!ht]
\begin{center}
\begin{tabular}{rrrrrrr} \hline
 N      & 20 & 40 & 80 & 160 & 320 & 640 \\ \hline 
$|k_1 -2 |$ & 1.07E-3  & 1.54E-4  & 4.53E-6  & 1.03E-6  &  2.12E-7   & 2.37E-8 \\  
$|k_2 -3 |$ & 1.99E-3  & 5.03E-4  & 1.09E-5  & 2.46E-5  &  5.16E-6   & 1.16E-6 \\ \hline
\end{tabular}
\caption{Numerical error of coefficients for two modes of low regular solutions in 1-D Poisson equation. 
Interface formulation is solved using second-order MIB method. $N$ is the number of uniform grid divisions.}
\label{tab:1d_poisson_error_k_2}
\end{center}
\end{table}

\begin{table}[!ht]
\begin{center}
\begin{tabular}{rrrrrrr} \hline
\multirow{2}{*}{N} &  \multicolumn{3}{c}{IF-MIB2} & \multicolumn{3}{c}{IF-MIB4} \\ \cline{2-4} \cline{5-7}
        &  $\| e \|_{\infty}$ & ROC & $|k-1|$ & $\| e \|_{\infty}$ & ROC & $|k-1|$   \\ \hline 
  20   &  7.65E-3  &      & 1.08E-2 & 6.59E-4  &       & 4.55E-4 \\ 
  40   &  2.28E-3  & 1.75 & 3.23E-3 & 4.97E-5  & 3.75  & 8.27E-5  \\ 
  80   &  6.38E-4  & 1.84 & 9.03E-4 & 3.38E-6  & 3.88  & 6.49E-6 \\ 
  160  &  1.53E-4  & 2.06 & 2.17E-4 & 2.11E-7  & 4.00  & 2.55E-7  \\ 
  320  &  3.10E-5  & 2.30 & 4.38E-5 & 1.43E-8  & 3.88  & 5.38E-8  \\  
  640  &  7.76E-6  & 2.00 & 3.48E-6 & 9.56E-10  & 3.90 & 1.86E-9   \\  \hline
\end{tabular}
\caption{Numerical error of ${u}$, its rate of convergence (ROC), and absolute error of $k$ for second-order 
and forth-order MIB methods with the interface formulation (IF-MIB2, IF-MIB4) on the 2-D Poisson equation 
with $u_s(x)=(x^2+y^2)^{1/4}$. $N$ is the same number of uniform grid divisions in $x$ and $y$ directions.}
\label{tab:2d_poisson_error_u}
\end{center}
\end{table}

\begin{figure}[!ht]
\begin{center}
\includegraphics[width=6.5cm]{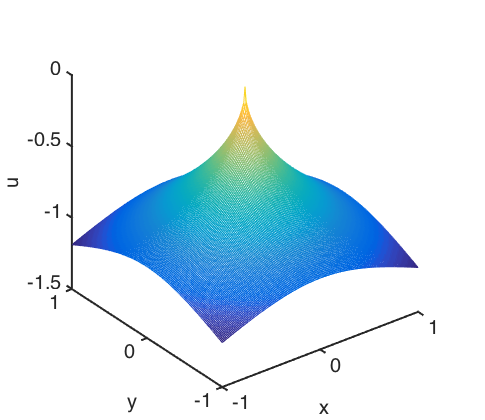}
\includegraphics[width=6.5cm]{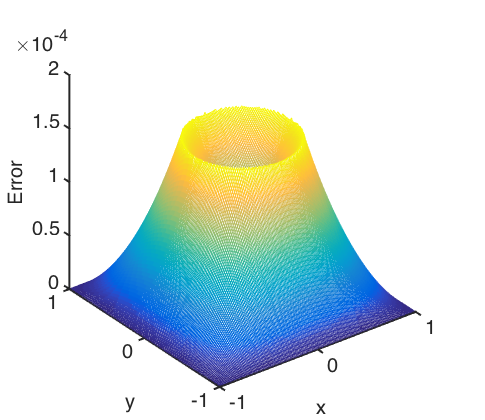}
\caption{Computational results on a $160 \times 160$ uniform mesh for 2-D Poisson equation with the analytical 
solution $u(x,y) = (x^2+y^2)^{1/4}$. The interface formulation and second order MIB method are used. Left: Negative 
of the numerical solution to feature the spike pointing upward; Right: Numerical error featuring the concentration 
of the error near the chosen interface $r=0.5$.}
\label{fig:error_2d_poisson}
\end{center}
\end{figure}

We proceed further to numerical experiments with a 2-D Poisson equation $\Delta u(x,y) = f(x,y)$ defined in the square domain $[-1,1] \times [-1,1]$.
The analytical solution is chosen to be $u(x,y) = (x^2 + y^2)^{1/4}$, which amounts to $r^{1/2}$ in the polar coordinate, 
featuring a low regularity at $(x,y) = 0$. We then choose a single mode $u_s(x,y) = (x^2 + y^2)^{1/4}$ and hence the analytical 
solution of $k$ is $1$. The interface is chosen to be a circle centered at the origin with a radius $r=0.5$. Inside the
circle a function $ku_s(x,y)$ is subtracted to generate the interface formulation, 
which is then solved using second and fourth order MIB methods. We collect the 
$L_{\infty}$ errors of ${u}(x)$ and the absolute error of $k$ in Table \ref{tab:2d_poisson_error_u}. 
These data evidence the comparable errors in $k$ and ${u}$ in the numerical solutions of 
the interface formulation. The left chart in Fig.(\ref{fig:error_2d_poisson}) highlights the low regularity of the solution at the 
origin, and the concentration of the numerical error near the chosen interface $r=0.5$ is shown in the right chart there. The consistent
second and forth order convergence observed from numerical solutions prove that our interface formulation can be utilized to remove
the low regular solution to generate an interface problem solvable by high order interface methods.

We shall finally make remarks on the construction of regular resolution modes $u_j^r(x)$ in $\Omega_1$ near $x_0$. When the Green's function 
$G$ of Eq.(\ref{eqn:model_problem}) is given, analytically or numerically, one can first specify a sequence representation of solution 
$\sum_{j=1}^{N_r} u^r_j(x)$ on the interface $\Gamma$. Then using a proper boundary integral with the given Green's function we can  
calculate $u^r_j(x)$ inside $\Omega_1$ and the normal derivative of $u^r_j(x)$ on $\Gamma$ (via Dirichlet-Neumann mapping). For 1-D
Poisson equation, this means we only need to specify $u^r(x)$ at two ends of $\Omega_1$. For 2-D or 3-D, one could use the spherical
harmonic functions as the basis to represent $u^r_j(x)$ on $\Gamma$. Direct expansion of $u^r_j(x)$ using Taylor polynomials with respect 
to $x_0$ is not recommended because the accuracy of the regular solution expansion is limited by the size of $\Omega_1$. On one hand,
we want to have $\Omega_1$ small to have highly accurate representation of regular solution components, but on the other hand, a smaller
$\Omega_1$ will increase the numerical errors in $\Omega_2$ because the truncation error there gets larger when $\Gamma$ is closer to $x_0$.

\section{Conclusions}
Many PDEs have solutions of low regularity whose fundamental modes can be analytically calculated, usually
with idealised boundary conditions. These low regular solution are generally localized but will degrade 
the global convergence of numerical methods when PDEs are solved with real boundary conditions. 
Here in this paper, we set the solution near the localized region as linear combination of known modes
of low regular solutions. By subtracting this combination from the solution in the chosen region we
obtained an interface problem with unknown interface conditions. This interface problem can be discretized by
using general interface methods. We designed two least square problems for the unknown coefficients in 
the linear combination. The solved coefficients with provide approximated interface conditions for 
the PDE, which can then be solved with high order numerical methods. Numerical experiments on 1-D and 2-D Poisson 
equations verified the usefulness of the designed interface formulation, as seen in the consistent second and forth
order convergence with Matched Interface and Boundary (MIB) methods. The errors in the coefficients of linear 
combination of low regular solutions and in the solution of PDEs are consistent, with a similar rate of convergence.

Our method could be regarded as a special domain decomposition method. Indeed, one could approximate the solution 
in the chosen region as the linear combination of known modes of low regular solutions and to solve the PDE
in this domain using collocation method. Standard numerical methods can be used in the complement domain,
and the boundary conditions will be exchanged on the interface between these two domains. This will lead to 
iterations between the numerical solutions between the two subdomains. Our interface formulation, 
nevertheless, avoids the kind of iterations by transferring the unknown coefficients 
of the linear combination into the interface conditions and by establishing least-square problems for the coefficients. 
When the low regular solution is completely known, it is not necessary to solve the unknown coefficients any more
with the interface formation. This is the case in the regularized Poisson-Boltzmann equation as investigated
in \cite{ChernY2003a,GengW2007a,HolstM2012a}, where the Coulomb potential induced by individual singular charges are 
completely known. We expect the promising applications of this interface formulation to linear elastic 
fracture mechanics where the solutions have various of modes of low regular solutions at the crack tip. 


\section*{References}
\bibliography{crack}
\end{document}